\documentclass[11pt]{article}

\usepackage{amsfonts,amsmath,graphicx}
\graphicspath{{.}}
   \makeatletter
   \let\accent@spacefactor\relax
   \makeatother
   \usepackage[french]{babel}
\usepackage{amsfonts}
\usepackage{amssymb}
\usepackage{amsmath}   %multiline etc. Subsumes amstext,amsbsy,amsopn
\usepackage{amsthm}    %provide proof env
\usepackage{amscd}     %commutative diagrams
\usepackage{euscript}
\usepackage[latin1]{inputenc}
\usepackage[french]{babel}

\newtheorem{defi}{D\'{e}finition}[section]
\newtheorem{pro}[defi]{Proposition}

\newtheorem{lem}[defi]{Lemme}
\newtheorem{theo}[defi]{Th\'{e}or\`{e}me}

\newtheorem{coro}[defi]{Corollaire}
\newtheorem{exa}[defi]{Exemple}

\newtheorem{rem}[defi]{Remarque}

\title{
{\bf Caractères numériques.}}
\author{Mireille Martin-Deschamps}
\date{janvier 2004}

\begin{document}
\maketitle

\begin{abstract}
La postulation des sous-schémas Arithmétiquement Cohen-Macaulay (ACM) de codimension 2 de  l'espace projectif ${\mathbb P}^N_k$ est bien connue, et a donné lieu à différentes approches : caractère numérique de Gruson/Peskine, $h$-vecteur, caractère de postulation de Martin-Deschamps/Perrin... Le premier but de cet article est d'établir l'équivalence de ces notions.

Le deuxième but, et le plus important, est d'étudier la postulation des sous-schémas ACM de  codimension 3 de ${\mathbb P}^N$. Pour cela on utilise la description due à Macaulay des fonctions de Hilbert des algèbres quotient d'un anneau de polynômes. On donne, par itération sur le nombre de variables, une nouvelle interprétation de la croissance de ces fonctions.

\medskip

The postulation of Arithmetically Cohen-Macaulay (ACM) subschemes of the projective space  ${\mathbb P}^N_k$ is well-known in the case of codimension 2. There are many different ways of recording this numerical information : numerical character of Gruson/Peskine, $h$-vector, postulation character of Martin-Deschamps/Perrin... The first aim of this paper is to show the equivalence between these notions.

The second, and most important aim, is to study the postulation of codimension 3 ACM subschemes of 
${\mathbb P}^N$. We use a result of Macaulay which describes all the Hilbert functions of the quotients of a polynomial ring. By iterating the number of variables, we obtain a new form of the growth of these functions.

\end{abstract}
\setcounter{section}{-1}
\section{Introduction}
\label{intro}

Soit $k$ un corps et  ${\mathbb P}^N_k$ l'espace projectif de dimension $n$ sur $K$. Pour classifier les sous-schémas de  ${\mathbb P}^N_k$ on leur associe des invariants numériques. Parmi ces invariants un des plus classiques est la postulation, qui donne pour chaque degré $d$ le nombre d'hypersurfaces indépendantes de degré $d$ contenant le sous-schéma considéré.

Le calcul de la postulation est un problème dont la complexité croit avec l codimension du sous-schéma.
 
 En codimension 1, la postulation d'une hypersurface $X$ est entièrement déterminée par son degré $d$ :  
puisque  le faisceau d'idéaux ${\mathcal I}_{X}$ qui la définit est  isomorphe à ${\mathcal O}_{\mathbb P}(-d)$, on a pour tout $n$, 
$h^0{\mathcal I}_{X}(n)=h^0 {\mathcal O}_{\mathbb P}(n-d)={{n-d+N}\choose{N}}$.

En codimension 2, et pour des sous-schémas localement Cohen-Macaulay, diverses notions ont  été introduites pour décrire cette postulation, le type numérique d'Ellingsrud \cite {E}, 
le caractère de postulation de \cite {MDP}, le caractère numérique de \cite {GP}, le h-vecteur (\cite {Mi}), ces deux dernières notions n'étant définies que pour des sous-schémas Arithmétiquement Cohen-Macaulay (ACM). Dans le cas des sous-schémas ACM, ces notions sont bien évidemment équivalentes.

En codimension au moins 3, aucun résultat n'est connu. Calculer la postulation de $X$ revient à calculer la fonction de Hilbert de l'algèbre graduée $k[X_0,\dots,X_N]/I_X$, où $I_X$ est l'idéal homogène saturé de $X$. On dispose d'un résultat très général du à Macaulay \cite {Mac} qui décrit toutes les fonctions de Hilbert des algèbres graduées quotient d'un anneau de polynômes, mais ce résultat, qui caractérise la ``croissance'' des fonctions considérées, n'est guère parlant. L'un des buts de cet article est de ``décrypter'' ce résultat et d'en donner une interprétation  qui soit plus utilisable dans la pratique. En particulier on pourra ainsi décrire toutes les postulations des sous-schémas ACM de codimension 3.

\medskip

Au premier paragraphe, on définit le caractère de postulation d'un sous-schéma fermé de l'espace projectif  ${\mathbb P}^N$ et sa variation par biliaison élémentaire Gorenstein.

Le deuxième paragraphe est consacré au résultat principal de cet article (\ref {theo principal}) qui donne une nouvelle caractérisation des fonctions de Macaulay.

Le troisième paragraphe applique ce résultat aux sous-schémas ACM de codimension 3. En particulier on en déduit au  quatrième paragraphe le calcul des degrés et genres des courbes ACM de ${\mathbb P}^4$ de degré inférieur ou égal à 10.

\subsection*{Notations}

On d\'esigne par $k$ un corps alg\'ebriquement clos, par ${\mathbb P}^N_k$ ou plus simplement
${\mathbb P}^N$ l'espace projectif de dimension $N\geq 2$ et par $S$
l'anneau de polyn\^omes
$k[X_0,\ldots,X_N]$.    Si ${\mathcal F}$ est un
 ${\mathcal O}_{\mathbb P}$-module on note $h^i {\mathcal F}$ la
dimension de l'espace vectoriel $H^i {\mathcal F}$.

Soient $X$ un sous-sch\'ema fermé de ${\mathbb P}^N$ et  ${\mathcal I}_X$ son 
faisceau d'id\'eaux, on d\'esigne par  $s_0(X)$ le plus petit degré d'une hypersurface contenant $X$, c'est-à-dire $s_0(X) = \inf \{\, n \in {\mathbb Z } \mid
h^0{\mathcal I}_X(n) \neq 0 \,\}$.

Soit $f : {\mathbb Z }\rightarrow {\mathbb Z}$ une application. On d\'efinit :

\begin{enumerate}
\item sa 
diff\'erence premi\`ere
$\partial f$ par $\partial f(n) = f(n) - f(n-1)$,
\item dans la cas où $f$ 
est nulle pour $n\ll 0$, sa primitive $f^{\sharp} $ par
 $f^{\sharp} (n) = \sum_ {k\le n} f(k) $,
\item sa borne supérieure si elle existe par $\sup f= \sup \{\, n \in {\mathbb Z } \mid
f(n) \neq 0 \,\}$,
\item sa décalée $f[d]$ par $f[d](n)=f(n+d)$.
\end{enumerate}

On rappelle qu'une fonction $f :{\mathbb Z }\rightarrow {\mathbb Z}$ , \`a support fini, et telle que l'on
ait $\sum _{n \in {\mathbb Z }} f(n) = 0$ est appel\'ee un {\bf
caract\`ere} \cite{MDP}, et que la diff\'erence premi\`ere d'une fonction \`a support fini est un caract\`ere.
On fait les conventions suivantes sur les  coefficients binômiaux :
\begin{equation*}
{{n}\choose{p}}=0  \; \textrm{pour } n\in {\mathbb Z }, \, p\geqslant 0\, \textrm{et } n<p
\quad ,  \quad {{n-1}\choose{-1}}=\begin{cases}1 
     & \textrm{pour } n=0\\
     0 & \textrm{sinon }\\
 \end{cases}
 \end{equation*}
 Avec cette convention, la formule de Pascal ${{n}\choose{p}}={{n-1}\choose{p}}+{{n-1}\choose{p-1}}$ est valable pour tous $n>p\geq 0$.
 
 On représentera parfois une fonction à support fini $f$ de $\mathbb{N}$ dans $\mathbb{N}$ par la suite de ses valeurs $(f(0),f(1),\cdots, f(\sup f))$, et la fonction ${{n-a-1}\choose{-1}}$ par $1_{[a]}$.

\section{Caract\`ere de postulation}

\begin{defi}
\rm Soient $X$ un sous-sch\'ema de pure dimension M de ${\mathbb P}^N$ et ${\mathcal I}_X$ son faisceau d'id\'eaux. On d\'efinit son
caract\`ere de postulation $\gamma_X$ par la formule :
$$\gamma_X(n)=\partial^{M+2}(h^0{\mathcal I}_{X}(n)-h^0 {\mathcal O}_{\mathbb P}(n)).$$

\end{defi}

\begin{rem}
\rm On notera que la fonction $h^0{\mathcal I}_{X}(n)-h^0 {\mathcal O}_{\mathbb P}(n)$ est l'oppos\'ee
de la fonction de Hilbert de $X$, dimension de l'image de la fl\`eche de
restriction ~:
$$H^0 {\mathcal O}_{\mathbb P}(n) \rightarrow H^0{\mathcal O}_{X}(n)$$
et qu'elle vaut $-h^0 {\mathcal O}_{X}(n)$ si $h^1{\mathcal I}_{X}(n)$ est nul, en particulier si $X$ est arithmétiquement Cohen-Macaulay.
\end{rem}

\begin{rem}
\rm Si $X$ est dégénéré, c'est-à-dire contenu dans un sous-schéma linéaire ${\mathbb P}^{N'}$de ${\mathbb P}^N$ avec $N'<N$, on voit facilement que son caractère de postulation est le même, calculé dans ${\mathbb P}^N$ ou dans ${\mathbb P}^{N'}$.
\end{rem}

\begin{exa}
\label{ex.hypersurface}
\rm Lorsque $X$ est une hypersurface de degré $d$ de ${\mathbb P}^N$ on a :

\begin{equation*}
\begin{split}
\gamma_X(n)&=\partial^{N+1}(h^0 {\mathcal O}_{\mathbb P}(n-d)-h^0 {\mathcal O}_{\mathbb P}(n))\\
&=\partial^{N+1}\left[{{n-d+N}\choose{N}}-{{n+N}\choose{N}}\right]\\
&={{n-d-1}\choose{-1}}-{{n-1}\choose{-1}}\\
\end{split}
\end{equation*}
donc la fonction $\gamma_X$ ne prend que deux valeurs non nulles, $\gamma_X(0)=-1$ et
$\gamma_X(d)=1$.

\end{exa}

\begin{pro}
\label{proprietes de gamma}
Soit $X$ un sous-sch\'ema de pure dimension M de ${\mathbb P}^N$ et $ \gamma= \gamma_X$ son caractère de postulation. Il a les propriétés suivantes :
\begin{enumerate}
\item [i)]$\sum _{n \in {\mathbb Z }} \gamma(n) = 0$ et donc $\gamma$ est un caract\`ere,
\item[ii)] $\gamma(n) = 0$ pour $n<0$,
\item [iii)]$\gamma(n) = -{{n+N-M-2} \choose {N-M-2}}$ pour $0 \leq n < s_0(X) = \inf \{\, n \in  {\mathbb Z } \mid h^0{\mathcal I}_X(n) \neq 0 \,\}$,
\item [iv)]$\gamma(s_0) > -{{s_0+N-M-2} \choose {N-M-2}}$.

\end{enumerate}

\end{pro}

\begin{proof}

 La fonction
$h^0{\mathcal I}_{X}(n)-h^0 {\mathcal O}_{\mathbb P}(n)$ est nulle pour $n<0$ et est égale à un polynôme de degré $M$
 pour $n \gg 0$. Sa diff\'erence $(M+1)$-ième est donc \`a support fini,
et $\gamma$ est un caract\`ere.

Pour $0 \leq n < s_0(X)$, on a :
$$\gamma(n)=-\partial^{M+2}h^0 {\mathcal O}_{\mathbb P}(n)
=-\partial^{M+2}{{n+N}\choose{N}}=-{{n+N-M-2} \choose {N-M-2}}$$
et
$$\gamma(s_0)=(h^0{\mathcal I}_{X}(s_0)-\partial^{M+2}h^0 {\mathcal O}_{\mathbb P}(s_0)
=h^0{\mathcal I}_{X}(s_0)-{{s_0+N-M-2} \choose {N-M-2}}$$
\end{proof}

\begin{pro}
Le caractère  $\gamma_X$ détermine la postulation et le polynôme de Hilbert de $X$.
\end{pro}

\begin{proof}
  La fonction $h^0{\mathcal I}_{X}(n)-h^0 {\mathcal O}_{\mathbb P}(n)$ est la primitive $(M+2)$-ième de $\gamma_X$, et on montre facilement  qu'on a \cite {MDP} :
$$h^0{\mathcal I}_{X}(n)-h^0 {\mathcal O}_{\mathbb P}(n)=\sum _{k \in {\mathbb Z }}{{n-k+M+1}\choose{M+1}}\gamma_X(k).$$
Soit $P_X$ le  polynôme de Hilbert de $X$. Pour $n\gg0$ on a : $$h^0{\mathcal I}_{X}(n)-h^0 {\mathcal O}_{\mathbb P}(n)= -h^0{\mathcal O}_{X}(n)=-P_X(n)$$
et $${{n-k+M+1}\choose{M+1}}={{(n-k+M+1)(n-k+M)\cdots(n-k+1)}\over {(M+1)!}}.$$
On en déduit l'égalité des polynômes :
$$P_X(n)=-\sum _{k \in {\mathbb Z }}{{(n-k+M+1)(n-k+M)\cdots(n-k+1)}\over {(M+1)!}}\gamma_X(k).$$
  \end{proof}

\begin{coro}
Soit $X$ un sous-schéma de ${\mathbb P}^N$ de degré $d$. On a :
$$d=\sum _{k \in {\mathbb Z }}k\gamma_X(k).$$
Soit $C$ une
courbe localement Cohen Macaulay de degré
$d$ et genre
$g$ de 
${\mathbb P}^N$. On a :
$$ g-1=\sum _{k \in {\mathbb Z }}
{{(k-1)(k-2)\over 2}}\gamma_C(k).$$
Soit $S$ une surface lisse de degré $d$ et genre arithmétique $p_a$ de ${\mathbb P}^N$ et $\delta$ le degré de son diviseur canonique. On a :
$$\delta=\sum _{k \in {\mathbb Z} }(k^2-4k)\gamma_S(k) \qquad
1+p_a=\sum _{k \in {\mathbb Z }}{{(k-3)(k-2)(k-1)}\over {6}}\gamma_S(k).$$

\end{coro}

\begin{proof} 
  En identifiant les coefficients de $n^M$ dans les deux membres de l'égalité :
  $$P_X(n)=-\sum _{k \in {\mathbb Z }}{{(n-k+M+1)(n-k+M)\cdots(n-k+1)}\over {(M+1)!}}\gamma_X(k)$$
  et en tenant compte de $\sum _{k \in {\mathbb Z }}\gamma_X(k) =0$ on obtient :
  
  $${{d}\over {M!}}=-\sum _{k \in {\mathbb Z }}{{(M+1-k)+\cdots+(1-k)}\over {(M+1)!}}\gamma_X(k)=\sum _{k \in {\mathbb Z }}{{(M+1)k}\over {(M+1)!}}\gamma_X(k).$$
  
 On a  :
\begin{equation*}
\begin{split}
P_C(n)&=nd+1-g =-\sum _{k \in {\mathbb Z }}{{(n-k+2)(n-k+1)}\over {2}}\gamma_C(k)\\
& =-\sum _{k \in {\mathbb Z }}{{(n^2-(2k-3)n+(k-1)(k-2)}\over 2})\gamma_C(k).\\
\end{split}
\end{equation*}

De même on a :
\begin{equation*}
\begin{split}
P_S(n)&={1\over 2}n^2 d-{1\over 2}n\delta+1+p_a =-\sum _{k \in {\mathbb Z }}{{(n-k+3)(n-k+2)(n-k+1)}\over {6}}\gamma_S(k)\\
& =-\sum _{k \in {\mathbb Z }}{{(n^3-(3k-6)n^2+(3k^2-12 k+11)n-(k-3)(k-1)(k-2)}\over 6})\gamma_S(k).\\
\end{split}
\end{equation*}
En tenant compte de $\sum _{k \in {\mathbb Z }}\gamma_C(k) =0$  et $\sum _{k \in {\mathbb Z }}\gamma_S(k) =0$, on en déduit les égalités annoncées.
\end{proof}

 \begin{pro}
 \label{definition de s1}
 Soient $X$ un sous-sch\'ema de pure dimension M de ${\mathbb P}^N$ et $ \gamma_X$ son caractère de postulation. Soit $s_1$ le plus petit degré d'une hypersurface contenant $X$ et ne contenant pas d'hypersurface de degré $s_0(X)$  contenant $X$. Alors on a :
 $$s_1= \inf
\big\{\, n \geq s_0 \mid \gamma_C(n) >{{n-s_0+N-M-2}\choose{N-M-2}}-{{n+N-M-2}\choose{N-M-2}} \,\big\}.$$

 \end{pro}
\begin{proof}
Soit $s_0=s_0(C)$. Pour $n<s_1$ on a :
$$h^0{\mathcal I}_C(n) ={{n-s_0+N}\choose{N}}\quad \textrm{et} \quad h^0{\mathcal I}_C(s_1) >{{s_1-s_0+N}\choose{N}}.$$
On en déduit, pour $s_0\leq n<s_1$, 
$$\gamma_X(n)={{n-s_0+N-M-2}\choose{N-M-2}}-{{n+N-M-2}\choose{N-M-2}} .$$

On a aussi : 
\begin{equation*}
\begin{split}
\gamma_X(s_1)&=h^0{\mathcal I}_X(s_1) -{{s_1-s_0+N}\choose{N}}+{{n-s_0+N-M-2}\choose{N-M-2}}-{{n+N-M-2}\choose{N-M-2}}\\
&>{{n-s_0+N-M-2}\choose{N-M-2}}-{{n+N-M-2}\choose{N-M-2}}.
\end{split}
\end{equation*}

\end{proof}
 
 \begin{rem}
\rm Dans le cas de la codimension 2, on obtient :
$$s_1= \inf
\{\, n \geq s_0 \mid \gamma_X(n) >0\,\}$$ et dans le cas de la codimension 3 :
$$s_1= \inf
\{\, n \geq s_0 \mid \gamma_X(n) >-s_0 \,\}.$$
\end{rem}

\subsection*{Lien avec les résolutions}

Soit $X$ un sous-sch\'ema de pure dimension M de ${\mathbb P}^N$ et $I_X$ son idéal saturé dans $S$. Puisque $I_X$ est de profondeur au moins 1 sur $S$, sa dimension projective est au plus $N-1$. Il existe donc une résolution graduée (qu'on peut choisir minimale) de $I_X$ par des $S$-modules libres gradués :

\[ 0\to L_{N-1}\to L_{N-2} \to \cdots \to L_0 \to I_X\to 0\]
Pour tout $i\in[1,N-1]$, on écrit $L_i=\oplus_{n\in {\mathbb Z }} {\mathcal O}_{\mathbb P}(-n)^{l_i(n)}$ et on définit la fonction $r_X$ par 
$$r_X(n)=\sum_{i=0}^{N-1}(-1)^i l_i(n)-{{n-1}\choose{-1}}.$$
On remarque que pour des raisons de rang,  $r_X$ est un caractère.
\begin{pro}
On a : $$r_X=\partial^{N-M-1} \gamma_X \quad \quad  \gamma_X(n)=\sum _{k \in {\mathbb Z }}{{n-k+N-M-2}\choose{N-M-2}}r_X(k).$$
\end{pro}

\begin{proof}
Notons ${\mathcal L}_i$ le faisceau dissocié associé au $S$-module libres gradué $L_i$. On a :
$$h^0{\mathcal L}_{i}(n)=\sum _{k \in {\mathbb Z }}{{n-k+N}\choose{N}}l_i(k).$$
De la  résolution de $I_X$, on déduit :

\begin{equation*}
\begin{split}
h^0{\mathcal I}_{X}(n)&=\sum_{i=0}^{N-1}(-1)^i h^0{\mathcal L}_{i}(n)\\
&=\sum_{i=0}^{N-1}(-1)^i \sum _{k \in {\mathbb Z }}l_i(k){{n-k+N}\choose{N}}\\
&=\sum _{k \in {\mathbb Z }}\left[r_X(k)+{{k-1}\choose {-1}}\right]{{n-k+N}\choose{N}}\\
&=\sum _{k \in {\mathbb Z }}r_X(k){{n-k+N}\choose{N}}+{{n+N}\choose {N}}.\\
\end{split}
\end{equation*}

En différenciant $M+2$ fois l'égalité :
$$h^0{\mathcal I}_{X}(n)-h^0 {\mathcal O}_{\mathbb P}(n)=\sum _{k \in {\mathbb Z }}r_X(k){{n-k+N}\choose{N}}$$
on obtient :
$$ \gamma_X(n)=\sum _{k \in {\mathbb Z }}{{n-k+N-M-2}\choose{N-M-2}}r_X(k).$$
En différenciant encore $N-M-1$ fois, on obtient : $\partial^{N-M-1} \gamma_X=r_X$.
\end{proof}

\subsection*{Variation par biliaison élémentaire Gorenstein}

Rappelons la notion de biliaison \'el\'ementaire Gorenstein, introduite par Hartshorne. 

\begin{defi} 
\label{bil-elem}
\rm \cite{H1} Soient $X$ et $X'$ deux sous-sch\'emas fermés de pure dimension M sans composante immergée de ${\mathbb P}^N$ trac\'es sur un sous-schéma fermé ACM, $Y$, de dimension $M+1$ satisfaisant $G_1$ (Gorenstein en codimension 1), et soit $h\in {\mathbb Z}$. On dit que $X'$ est obtenu par une
biliaison \'el\'ementaire Gorenstein de hauteur $h$ sur $Y$  \`a partir de $X$ (ascendante si
$h>0$, descendante si $h<0$) si on a une \'equivalence lin\'eaire de diviseurs
g\'en\'eralis\'es $X' \sim X + hH$, c'est-à-dire un isomorphisme de faisceaux d'idéaux relatifs 
${\mathcal I}_{X/Y}\simeq {\mathcal I}_{X'/Y}(h)$.\\
\end{defi}

 Du point de vue cohomologique, si on a $h>0$, une 
biliaison de hauteur $h$ sur $Y$ 
est \'equivalente \`a $h$ biliaisons de hauteur 1 sur $Y$. On pourra donc  se borner
 \`a faire le calcul dans ce dernier cas.

\begin{pro}
 \label{variation Gorenstein}
Si $X'$ est obtenu par une
biliaison \'el\'ementaire Gorenstein de hauteur $1$ sur $Y$  \`a partir de $X$, on a :
$$\gamma_{X'}(n)=\gamma_X(n-1)+\gamma_Y(n).$$
\end{pro}
\begin{proof}

On utilise la relation ${\mathcal I}_{X/Y}\simeq {\mathcal I}_{X'/Y}(h)$
   et les deux suites exactes :
$$0 \to   {\mathcal I}_{Y} \to  {\mathcal I}_{X}\to{\mathcal I}_{X/Y} 
\to0 $$
$$0 \to {\mathcal I}_{Y} \to {\mathcal I}_{X'}\to {\mathcal I}_{X'/Y}
 \to 0 . $$

On en déduit  :
$$h^0 {\mathcal I}_{X'} (n) = h^0 {\mathcal I}_X  (n-h) + h^0 {\mathcal I}_{Y} (n) -h^0 {\mathcal I}_{Y} (n-h) $$
 et en différentiant $M+2$ fois :
$$\gamma_{X'} (n) = \gamma_{X} (n-h) + \gamma_{Y} ^{\sharp}(n)-\gamma_{Y} ^{\sharp}(n-h).$$
Pour $h = 1$ :
$$\gamma_{X'} (n) = \gamma_{X} (n-1) +\gamma_{Y} ^{\sharp}(n)-\gamma_{Y} ^{\sharp}(n-1)= \gamma_{X} (n-1) +\gamma_{Y} (n).$$
   
   \end{proof}

\begin{rem}
\rm Ce résultat est classique pour les sous-schémas de codimension 2 de ${\mathbb
P}^N$. Dans ce cas, $Y$ est une hypersurface de degré $s$ et la notion de biliaison élémentaire Gorenstein co\"\i ncide avec celle de biliaison élémentaire intersection complète. 
 Si $X'$ est obtenu par une
biliaison \'el\'ementaire de hauteur $h$ sur $Y$ \`a partir de $X$, on déduit de \ref{ex.hypersurface} qu'on a :
     \begin{equation*}
\gamma_{X'}(n)-\gamma_{X}(n-1)=\begin{cases}-1  & \textrm{pour } n=0\\
     1 & \textrm{pour } n=s\\
     0 & \textrm{sinon }\\
 \end{cases}
 \end{equation*}
 
\end{rem}

\section{Sous-schémas ACM et conditions  de Macaulay}

Dans ce paragraphe, nous allons étudier les sous-schémas  de l'espace projectif
${\mathbb P}^N$ qui sont arithm\'etiquement Cohen-Macaulay (en abrégé ACM),
c'est-à-dire dont l'anneau gradué
$S/{I_X} $ est Cohen-Macaulay (de dimension $1+\textrm{dim} X$).

\begin{defi} 
\rm Soit  $X$ un sous-sch\'ema ferm\'e ACM de ${\mathbb
P}^N$ de dimension $M$ et $I_X$ son idéal saturé. 
Si $L_1,\ldots,L_{M+1}$
 sont des formes linéaires générales, elles forment une suite régulière pour $S/{I_X} $. On en déduit qu'on a des suites exactes :
$$0\to S/{I_X}(-1)\to   S/{I_X}\to S/{I_X+(L_1)}$$
$$0\to S/{I_X+(L_1,\dots,L_i)}(-1)\to   S/{I_X}\to S/{I_X+(L_1,\dots,L_{i+1})}\to 0$$
de $S$-modules gradués, qui permettent de calculer la fonction de Hilbert $h_X$ de 
l'anneau artinien $S/{I_X+} (L_1,\ldots,L_{M+1})$. Cette fonction de Hilbert $h_X$ est le
{\bf h-vecteur}  de $X$. 

\end{defi}

\begin{rem}
\rm On voit facilement qu'on a :
$$h_X(n)=\partial^{M+1}(h^0 {\mathcal O}_{\mathbb P}(n)-h^0{\mathcal I}_{X}(n))$$
donc $\gamma_X=-\partial h_X$ et la  notion de h-vecteur pour un sous-schéma ACM est
équivalente à celle de caractère de postulation. De plus, si $d$ est le degré de $X$, on a l'égalité $d=\sum _{k \in {\mathbb Z }}k\gamma_X(k)=\sum _{k \in {\mathbb Z }}h_X(k)$.
\end{rem}

\begin{rem}
\rm Dans le cas des sous-schémas ACM de codimension 2, il existe une autre notion
\'equivalente, c'est le caract\`ere numérique de Gruson-Peskine \cite{GP}.
\end{rem}

L'anneau gradué $S/{I_X} $ d'un sous-schéma est une $k$-alg\`ebre de type fini,
engendr\'ee par sa composante de degr\'e 1. C'est ce que certains auteurs appellent
une {\bf G-alg\`ebre standard}. Les fonctions de Hilbert de ces algèbres ont été
caractérisées par Macaulay. 
 C'est pourquoi
toutes les variantes de la postulation des sous-schémas  ACM (caract\`ere de
Gruson-Peskine, caract\`ere de postulation, $h$-vecteur) peuvent
s'obtenir \`a partir des fonctions de Macaulay que nous définissons ci-dessous.
On se reportera pour les d\'etails \`a \cite{Mac} ou \cite{S}, qui reprend les
r\'esultats sous une forme plus moderne (voir aussi \cite{GMR}).

\subsection* {Fonctions de Macaulay}

\begin{pro}
\label{existence developpement i-binomial} Soient $\alpha$ et $i$ des entiers
strictement positifs. Alors
$\alpha$ peut s'\'ecrire de mani\`ere unique sous la forme:
$$\alpha={{m_i}\choose {i}} +{{m_{i-1}}\choose {i-1}}+\cdots +{{m_j}\choose
{j}}$$ avec $m_i>m_{i-1}>\cdots>m_j\geq j\geq 1$. 
\end {pro}

\begin{proof}
L'existence et l'unicit\'e d\'ecoulent des in\'egalit\'es suivantes, qui
caract\'erisent $m_i$) :
$${{m_i}\choose {i}}\leq \alpha<{{m_i+1}\choose {i}}={{m_i}\choose {i}} +
{{m_i-1}\choose {i-1}}+\cdots +{{m_i-i+1}\choose
{1}}+1.$$
\end{proof}

\begin{defi}
\label{developpement i-binomial}
\rm L'expression de $\alpha$ en fonction des coefficients binômiaux établie en
\ref{existence developpement i-binomial} est appelée le {\bf d\'eveloppement
i-bin\^omial} de $\alpha$.\\
On d\'efinit alors  :
$$\alpha^{<i>}={{m_i+1}\choose {i+1}} +{{m_{i-1}+1}\choose {i}}+\cdots
+{{m_j+1}\choose {j+1}}$$ 
(qui est le d\'eveloppement $(i+1)$-bin\^omial de $\alpha^{<i>}$) et $0^{<i>}=0$.
\end{defi}

\begin{exa}
 $25={6\choose 3}+{3\choose
2}+{2\choose 1}$, $\;25^{<3>}={7\choose 4}+{4\choose
3}+{3\choose 2}=41$.
\end{exa}

\begin{rem}
\rm
 Si $0<\alpha\leq i$, les coefficients bin\^omiaux qui
apparaissent dans le d\'eveloppement i-bin\^omial de $\alpha$ sont \'egaux
\`a 1. On en d\'eduit facilement qu'on a alors $\alpha^{<i>}=\alpha$.
\end{rem}

\begin{defi}
\label{conditions de croissance}
\rm Une fonction $h$ de ${\mathbb N}$ dans ${\mathbb N}$ satisfait aux conditions de croissance de Macaulay si elle vérifie
$h(0)=1$
et
$h(i+1)\leq h(i)^{<i>}$ pour tout $i\geq 1$. Si $a=h(1)$, on dira aussi que c'est une fonction de
Macaulay de type $a$. 
\end{defi}

Le r\'esultat de Macaulay est le suivant :

\begin{theo}
\label{theo Macaulay}
Une fonction satisfait aux conditions de croissance de Macaulay si et seulement si c'est la fonction de Hilbert d'une G-alg\`ebre standard.
\end{theo}

\begin{rem}
\rm Toute fonction $h$ binômiale, c'est-à-dire de la forme $h(n)={{p+n}\choose {n}}$, où
$p$ est un entier positif, est une fonction de Macaulay de type $p+1$. Plus
généralement, si une fonction $h$ est binômiale sur un intervalle $[n_0,n_1]$ (resp.
$[0,n_1]$), elle satisfait aux conditions de croissance sur l'intervalle
$[n_0+1,n_1]$ (resp. $[0,n_1]$).
\end {rem}

Nous allons mettre en \'evidence quelques propri\'et\'es simples des
fonctions de Macaulay.

\begin{pro}.  
\label{proprietes}
\begin{itemize}
\item[a)] Soient $\alpha$, $\beta$ et $i$ des
entiers strictement positifs avec $\alpha<\beta$. Alors on a
$\alpha^{<i>}<\beta^{<i>}$.\\
 \item[b)] Soit $h$ une fonction de Macaulay de type $a$. Alors
$h(n+1)\leq {a+n\choose {n+1}}$ pour tout
$n\geq 0$.\\
\item[c)] S'il existe $n>1$ tel que $h(n)=0$, alors $h(m)=0$ pour tout $m\geq
n$.\\
 \item[d)] S'il existe $n>1$ tel que $h(n)=1$, alors $h(m)\leq 1$ pour tout
$m\geq n$.\\
\item[e)]  S'il existe $n>1$ tel que $h(n)<n+1$, alors la fonction $h$
est d\'ecroissante pour 
$m>n$.
\end{itemize}
\end{pro}

\begin{proof}
 a)  Ecrivons les
d\'eveloppements i-bin\^omiaux :$$\alpha={{m_i}\choose {i}}
+{{m_{i-1}}\choose {i-1}}+\cdots +{{m_j}\choose {j}},\quad
\beta={{n_i}\choose {i}} +{{n_{i-1}}\choose {i-1}}+\cdots
+{{n_{j'}}\choose {j'}}.$$
Quitte \`a retrancher un m\^eme nombre \`a $\alpha$ et $\beta$, on peut
supposer qu'on a
$m_i\neq n_i$.

Si $m_i<n_i$, alors $\alpha<{{m_i+1}\choose {i}}\leq {{n_i}\choose
{i}}\leq \beta$. De m\^eme si  $m_i>n_i$, alors $\alpha>\beta$.
Donc si  $\alpha<\beta$, on a $m_i<n_i$, $m_i+1<n_i+1$ et
$\alpha^{<i>}<\beta^{<i>}$.

b) r\'esulte de a).

d) r\'esulte de $1^{<n>}=1$.

e) Si $h(n)<{n+1\choose n}$, alors $h(n)^{<n>}=h(n)$.
\end{proof}

\begin{rem}
\rm Soit $h$ une fonction de Macaulay de type $a\neq 0$, vu \ref {proprietes} b, il est naturel de lui associer
l'entier (éventuellement infini)
$s_0(h)= 
\inf \{\, n \in {\mathbb Z } \mid
h(n) <{{a+n-1}\choose {n}} \,\}$. Par définition il est  $>1$. S'il est infini, $h$ est binômiale : $h(n)={{a+n-1}\choose {n}} $.
\end{rem}

D'après \ref{theo Macaulay}, toute fonction de Macaulay à support fini est la
fonction de Hilbert d'une G-alg\`ebre standard de longueur finie. On montre en fait
que c'est le quotient d'un anneau de polynômes par un idéal monomial.
En rajoutant des variables on peut relever un tel idéal en un idéal saturé de
$k[X_0,\ldots,X_N]$ de la même codimension. On a donc montré le résultat suivant :

\begin {theo}
%%% Etudier le résultat de {Migliore-Nagel}
\label{Migliore-Nagel}
Tout fonction de Macaulay de type $a$ à support fini est le h-vecteur d'un sous-schéma ACM de
${\mathbb P}^N$ de codimension $a$.
\end{theo}

\begin{rem}
\label{borne inf}
\rm Si $X$ n'est pas dégénéré, le type de $h_X$ est égal à la codimension $a$ de $X$ et $s_0(h_X)=s_0(X)$. On peut alors donner une borne inférieure pour le degré de $X$ :
 $$d=\sum _{k \in {\mathbb Z }}h_X(k)\geqslant \sum _{0\leqslant k <s}h_X(k)=\sum _{0\leqslant k <s}{{a+k-1}\choose {k}}= {{a+s-1}\choose {s-1}}.$$
\end{rem}

Le lemme suivant établit un résultat technique qui sera utile dans la preuve du théorème principal  
\ref{theo principal}.

\begin{lem}
\label{lemme technique}
Soient $\alpha$, $\beta$ et $i$ des
entiers strictement positifs et 
$$\alpha={{m_i}\choose {i}} +{{m_{i-1}}\choose {i-1}}+\ldots +{{m_j}\choose
{j}}$$
le développement i-binomial de $\alpha$. Supposons qu'on ait $\beta<{{m_j}\choose
{j-1}}$. Alors on a :
$$(\alpha +\beta)^{<i>}=\alpha^{<i>}+\beta^{<j-1>}.$$
\end{lem}

\begin{proof}
Ecrivons le développement (j-1)-binomial de $\beta$ :
$$\beta={{m'_{j-1}}\choose {j-1}} +{{m'_{j-2}}\choose {j-2}}+\ldots +{{m'_k}\choose
{k}}$$
On a $\beta<{{m_j}\choose
{j-1}}$ donc $m'_{j-1}<m_j$ et l'écriture :
$$\alpha+\beta={{m_i}\choose {i}} +\ldots +{{m_j}\choose{j}}+
{{m'_{j-1}}\choose {j-1}}+\ldots +{{m'_k}\choose
{k}}$$
est  le développement i-binomial de $\alpha+\beta$.

Alors on a 
\begin{equation*}
\begin{split}
(\alpha +\beta)^{<i>}&={{m_{i+1}}\choose {i+1}} +\cdots
+{{m_{j+1}}\choose{j+1}}+{{m'_{j}}\choose {j}} +\cdots
+{{m'_{k+1}}\choose {k+1}}\\
&=\alpha^{<i>}+\beta^{<j-1>}\\
\end{split}
\end{equation*}
\end{proof}

Il est facile de décrire les fonctions de Macaulay de type 1 ou 2.

\begin{pro}
Une fonction $h$ est une fonction de Macaulay de type 1 si et seulement si elle est 
d\'ecroissante et prend seulement les valeurs 1 et 0. \\
Une fonction $h$ est une fonction de Macaulay de type 2 si et seulement si il existe un entier $s_0>
1$ (éventuellement infini) tel que
$h(n)=n+1$ pour $n< s_0$ et tel que $h$ soit  d\'ecroissante pour $n\geq
s_0$. On a alors $s_0=s_0(h)$.
\end{pro}
\begin{proof}
La première assertion est une conséquence immédiate de \ref{proprietes}.

Pour le type 2, on pose $s_0=1+  \sup \{ n \mid h(n)=n+1 \; \}$. Pour 
$n\geq s_0$,  le r\'esultat est encore une conséquence  de \ref{proprietes}.

Les réciproques sont immédiates.
\end{proof}

\rm  On retrouve ainsi immédiatement la positivité  du caractère d'un sous-sch\'ema ACM de  pure codimension 2  (cf. \cite {GP}, \cite{MDP}). Posons tout d'abord la définition suivante :

\begin{defi}
\label{definition caractere positif}
\rm Un caractère $\gamma$ est positif s'il vérifie $\gamma(n)=0$ pour $n<0$, $\gamma(0)=-1$, $\gamma(n)\geq 0$ pour tout $n\geq s_0(\gamma)=\inf \{\, n \in {\mathbb Z} \mid \gamma(n) \neq -1 \,\}$.
\end{defi}

\begin{rem}
\label{caractere positif-fonction de Macaulay}
\rm Soit $h$ une
fonction de Macaulay de type 2. D'après ce qui précède sa différence première $\partial h$
vérifie $\partial h(n)=1$ pour $0\leq n<s_0$ et $\partial h(n)\leq 0$ pour $n\geq
s_0$. 

On en déduit qu'une fonction $h$ à support fini est une fonction de Macaulay de type 2 si et seulement si l'opposée de sa différence première $-\partial h$
est un caractère positif $\gamma$ vérifiant $s_0(\gamma)\geqslant 2$. Dans ce cas, on a $s_0(\gamma)=s_0(h)$.

Si $h$ est  de type 1,  $-\partial h$  est le caractère (positif) de la section hyperplane d'une hypersurface de degré $d$. On a alors $s_0(\gamma)=1$ et $s_0(h)=d$. 

Si $h$ est  de type 0,  $-\partial h$  est le caractère (positif) de l'intersection de deux hyperplans. On a alors $s_0(\gamma)=1$ et $s_0(h)$ n'est pas défini. 

\end{rem}

\begin{theo}
   \label{caractere positif}
\cite {GP}  Soit $X$  un sous-sch\'ema ACM de  pure codimension 2 de ${\mathbb P}^N$. Alors son caractère $\gamma_X$ est positif. 
Inversement, soit $\gamma$ un caractère positif. Alors il existe un sous-sch\'ema X  de  pure codimension 2 et ACM de ${\mathbb P}^N$ tel qu'on ait  $\gamma=\gamma_X$.
\end{theo}

Dès qu'on a $h(1)\geq 3$, il est beaucoup plus difficile de décrire les fonctions de Macaulay. C'est ce que  nous tentons de faire avec le résultat suivant :

\begin{theo}
\label{theo principal}
Soit $h$ une fonction non binômiale vérifiant $h(0)=1$ et $h(1)=a>1$. Alors $h$ 
est une fonction de Macaulay si et seulement si il existe $r+1$ fonctions de
Macaulay ($r\geq 1$) $h_0, \ldots,h_r$ avec :
\begin{itemize}
\item[i)] $h_i(1)=a-1$ pour $0\leq i<r$, $h_r(1)\leq a-1$,
\item[ii)]  pour $1\leqslant i\leqslant r$, $\sup h_i <s_0(h_{i-1})-1$ ou $s_0(h_{i-1})$ est infini,
\item[iii)] $h=h_0+h_1[-1]+\cdots +h_r[-r]$.
\end{itemize}
De plus, $h_0, \ldots,h_r$ sont déterminés de manière unique par ces conditions et
on a $s_0(h)=r+1$.
\end{theo}

 La démonstration  va se faire en plusieurs étapes. On commence par un lemme :
 \begin{lem}
\label{proprietes de h}

Soient  $h_0, \ldots,h_r$ ($r\geq 1$) des fonctions de Macaulay vérifiant les
conditions de l'énoncé de \ref{theo principal}. Posons $t=\sup \{i\in \mathbb{N}\mid s_0(h_i) =+\infty \}$ si cet ensemble est non vide, $t=-1$ sinon. Alors $h_i$ est binômiale pour $0\leqslant i\leqslant t$.  De plus, la fonction   $h=h_0+h_1[-1]+\cdots +h_r[-r]$
est une fonction de Macaulay de type $a$ qui vérifie les propriétés suivantes :
\begin{itemize}
\item $s_0(h)=r+1$,
\item  $s_0(h_{t+1})=\inf \{\, n \in
{\mathbb N }
\mid h(n)<\sum_{i=0}^{t+1} {{a+n-2-i}\choose{n-i}} \,\}$,
\item $h(n)=\sum_{i=0}^{j} {{a+n-2-i}\choose{n-i}} \,\} +h'(n-j-1)$ où $h'$ est une fonction de Macaulay de type $\leqslant a$.
\end{itemize}
\end{lem}

\begin{proof}
On a $h(0)=h_0(0)=1$, $h(1)=h_0(1)+h_1(0)=a$.

Pour $1\leqslant i\leqslant r$, si $s_0(h_{i})$ est infini, ce qui revient à dire que $h_i$ est binômiale, $\sup h_i $ est  infini, donc d'après la condition ii) de \ref{theo principal}, $s_0(h_{i-1})$ est également infini. L'ensemble des indices $i$ tels que $s_0(h_{i})$ soit infini est donc, s'il est non vide, un intervalle $[0,t]$.

Pour tout $t+2\leqslant i<r$ on a $\sup h_i\geq s_0(h_{i})-1$, donc la condition  $\sup h_i
<s_0(h_{i-1})-1$ entraine que la suite des $s_0(h_{i})$ pour $t+1\leqslant i<r$ est strictement décroissante.
On a :
$$r+1\leq \sup h_r+r+1\leq s_0(h_{r-1})+r-1\leq \cdots \leq s_0(h_{t+1}).$$

Pour $n\leq r$, on a $n<s_0(h_{i})+i$ pour $i<r$ donc d'après la définition de
$s_0(h_{i})$ et en utilisant le fait que $h_i$ est une fonction de Macaulay de type
$a-1$ on a 
$h_i(n-i)={{a+n-i-2}\choose {n}}$ pour $i<r$. On a aussi :

\begin{equation*}
\begin{split}
h(n)&=h_0(n)+h_1(n-1)+\cdots +h_n(0)\\
&={{a+n-2}\choose {n}}+{{a+n-3}\choose {n-1}}+\cdots +1\\
&={{a+n-1}\choose {n}}\\
\end{split}
\end{equation*}

et :

\begin{equation*}
\begin{split}
h(r+1)&=h_0(r+1)+h_1(r)+\cdots +h_r(1)\\
&={{a+r-1}\choose {r+1}}+\cdots +{{a}\choose {2}}+h_r(1)\\
&\leq {{a+r-1}\choose {r+1}}+\cdots +{{a}\choose {2}}+{{a-1}\choose {1}}\\
&<{{a+r}\choose {r+1}}\\
\end{split}
\end{equation*}
donc on a $r+1=
\inf \{\, n \in {\mathbb Z } \mid
h(n) <{{a+n-1}\choose {n}} \,\}$, et la condition de croissance $h(n+1)\leqslant h(n)^{<n>}$ est vérifiée pour $0\leqslant n\leqslant r$.

Pour $n<s_0(h_{t+1})$ on a :
$$h(n)\geq h_0(n)+\cdots +h_{t+1}(n-t-1)=\sum_{i=0}^{t+1} {{a+n-2-i}\choose{n-i}}.$$

Pour $n\geq s_0(h_{t+1})$ on a $h_i(n)=0$ pour tout $i>t+1$, donc 
$$h(n)=
h_0(n)+\cdots +h_{t+1}(n-t-1)<\sum_{i=0}^{t+1} {{a+n-2-i}\choose{n-i}}.$$
 On en déduit d'une part qu'on a  :
 $$s_0(h_{t+1})=\inf \big\{\, n \in
{\mathbb N }
\mid h(n)<\sum_{i=0}^{t+1} {{a+n-2-i}\choose{n-i}} \,\big\},$$
  d'autre part que la condition de croissance est vérifiée pour $n\geq s_0(h_{t+1})$.

Pour $1\leqslant i \leqslant r-1$ et $s_0(h_{i})+i\leq n <s_0(h_{i-1})+i-1$, on a :
$$h(n)={{a+n-2}\choose{n}}+\cdots
+{{a+n-i-1}\choose{n-i+1}}+h_i(n-i)$$
et de plus $h_i(n-i)<{{a+n-i-1}\choose{n-i}}$.
D'après le lemme \ref{lemme technique} on en déduit :
 
$$
h(n)^{<n>}={{a+n-1}\choose{n+1}}+\cdots
+{{a+n-i}\choose{n-i+2}}+h_i(n-i)^{<n-i>}$$

De plus, pour tout $j\leqslant i-1$, on a $n<s_0(h_j)+1$, et on en déduit :
\begin{equation*}
\begin{split}
h(n+1)& =h_0(n+1)+\cdots +h_{i-1}(n+2-i)+h_i(n+1-i)\\
& \leqslant {{a+n-1}\choose{n+1}}+\cdots
+{{a+n-i}\choose{n-i+2}}+h_i(n-i)^{<n-i>}=h(n)^{<n>}.
\end{split}
\end{equation*}

Pour $r+1\leq n< s_0(h_{r-1})+r-1$ on a :
$$h(n)={{a+n-2}\choose{n}}+\cdots
+{{a+n-r-1}\choose{n-r+1}}+h_r(n-r)$$
et de plus $h_r(n-r)<{{a+n-r-1}\choose{n-r}}$ puisque $h_r$ est de type $a'\leqslant a-1$.
On conclut comme dans le cas précédent.

Posons $h'=h_{t+1}+\cdots +h_r[t+1-r]$ si $t\leqslant r-1$ et $h'=0$ si $t=r$,  c'est-à-dire qu'on a :
$$h(n)=\sum_{i=0}^{t} {{a+n-2-i}\choose{n-i}} \,\} +h'(n-t-1).$$
 Alors le résultat qu'on vient de montrer s'applique à $h'$ qui est une fonction de Macaulay, de type $a$ si $t+1<r$, de type $\leqslant a-1$ si $t+1=r$.
\end{proof}
\begin{proof} {\em de \ref {theo principal}}.

\subsubsection*{unicité.}

Supposons qu'on ait deux écritures :  
$$h=h_0+h_1[-1]+\cdots +h_r[-r]=
g_0+g_1[-1]+\cdots +g_{r'}[-r']$$
avec des fonctions vérifiant les conditions de l'énoncé. 

 On peut alors écrire, en regroupant les composantes binômiales : 
   $$h(n)=\sum_{i=0}^t {{a+n-i-2}\choose{n-i}}+h'(n-t-1)=\sum_{i=0}^{t'} {{a+n-i-2}\choose{n-i}}+g'(n-t'-1)$$
   où $h'$ et $g'$ sont des fonctions de Macaulay de type $\leqslant a$. 
   
   Si $t\neq -1$, $h(n)$ est équivalent à  $tn^{a-2}/(a-2)!$ quand $n$ tend vers $+\infty$. Si $t=-1$, $h(n)$ est d'ordre $<a-2$ en $n$ quand $n$ tend vers $+\infty$. On en déduit que $t=t'$, $h'=g'$, et en retranchant la somme des composantes binômiales, on se ramène au cas où $t=t'=-1$.
   
   D'après \ref {proprietes de h} on a $s_0(h)=r+1$ donc $r$ est déterminé par $h$, et $r=r'$.

De même, on a vu dans la preuve de  \ref {proprietes de h} que $s_0(h_0)$ est déterminé par $h$
donc $s_0(h_0)=s_0(g_0)$ et pour $n>s_0(h_0)$, $h(n)=h_0(n)=g_0(n)$. Les fonctions
$h_0$ et $g_0$ sont donc égales.

Soit alors $i_0\leq r-2$ tel qu'on ait $h_i=g_i$ pour tout $i\leq i_0$. Posons :
$$h'[-i_0-1]= h-h_0-h_1[-1]-\cdots +h_{i_0}[-i_0]$$ On a 
$$h'=h_{i_0+1}+h_{i_0+2}[-1]+\cdots
+h_r[i_0+1-r]=g_{i_0+1}+g_{i_0+2}[-1]+\cdots
+g_r[i_0+1-r]$$ 
et pour les mêmes raisons que ci-dessus on a :
$$s_0(h_{i_0+1})=s_0(g_{i_0+1})=\inf \Big\{\, n \in
{\mathbb Z }
\mid h'(n)< {{a+n-2}\choose{n}} \,\Big\}$$
et pour $n\geq s_0(h_{i_0+1})$, $ h'(n)=h_{i_0+1}(n)=g_{i_0+1}(n)$, donc
$h_{i_0+1}=g_{i_0+1}$.

On a donc $h_i=g_i$ pour tout $i\leq r-1$ et donc $h_r=g_r$ et l'écriture est unique.

\subsubsection*{existence quand $h$ est à support fini}
Montrons dans ce cas l'existence, pour $a$ fixé, par récurrence sur
$h^{\sharp}(\infty)=\sum_{n\in {\mathbb N }}h(n)\geq 1+a$.

\begin{itemize}
\item 1er cas : $h^{\sharp}(\infty)=1+a$. 

Alors on a $h(n)=0$ pour tout $n\geq 2$.
On définit les deux fonctions de Macaulay $h_0$ et $h_1$ par :

 \begin{equation*}
 h_0(n)=\begin{cases} 1 & \textrm{pour } n=0\\
 a-1 & \textrm{pour } n=1\\
 0 & \textrm{pour tout }n\geq 2\\
 \end{cases}
 \quad , \quad 
 h_1(n)=\begin{cases} 1 & \textrm{pour } n=0\\
 0 & \textrm{pour tout }n\geq 1\\
 \end{cases}
 \end{equation*}

  On a $s_0(h_{0})=2$, $\sup h_1=0$ et  $h=h_0+h_1[-1]$
donc les propriétés cherchées sont vérifiées.

\item  2ème cas : $h^{\sharp}(\infty)>1+a$. 
 
 Soit $N=\inf \{\, n \in
{\mathbb N }
\mid h(n)< {{a+n-2}\choose{n}} \,\}$.
Définissons les fonctions $h_0$ et $h'$ par :
 
 \begin{equation*}
h_0(n) =\begin{cases} {{a+n-2}\choose{n}} & \textrm{pour } n<N\\
h(n) & \textrm{pour }n\geq N\\
 \end{cases}\\
\quad , \quad h' =(h-h_0)[1].
 \end{equation*}

Par construction on a $s_0(h_{0})=N$.
Pour montrer que $h_0$ est une fonction de Macaulay (de type $a-1$) il suffit
de vérifier les conditions de croissance pour $n=N-1$, c'est-à-dire 
$h_0(N)\leq h_0(N-1)^{<N-1>}$. Or on a 
 $$h_0(N-1)={{a+N-3}\choose{N-1}}\quad 
 \textrm  {donc} \quad h_0(N-1)^{<N-1>}={{a+N-2}\choose{N}}>h_0(N).$$

Nous allons montrer maintenant que $h'$ est une fonction de Macaulay. C'est par
construction  une fonction de ${\mathbb N }$ dans $ {\mathbb N }$, on a :
$$h'(0)=h(1)-h_0(1)=1\quad , \quad h'(1)=h(2)-h_0(2).$$

\begin{itemize}
 \item Si $N=2$, $h(2)=h_0(2)$ et $h'(1)=0$.
\item Si $N>2$, 
 $$h_0(2)={a\choose 2}\quad \textrm  {et} \quad  {{a}\choose 2}\leq h(2)\leq {{a+1}\choose 2}\quad \textrm  {donc} \quad 0\leq h'(1)
\leq a.$$
 \end{itemize}

Puisque $\sup h'\leq N-2$ il suffit
de vérifier les conditions de croissance pour $n<N-2$.

Pour $n< s_0(h)$ : $$h'(n)= {{a+n}\choose {n+1}}- {{a+n-1}\choose
{n+1}}={{a+n-1}\choose {n}}$$
  est binômiale et les conditions de croissance sont
vérifiées. 

Pour $s_0(h)-1\leq n <N-2$, on a : 
  $$ h(n+1)< {{a+n}\choose
{n+1}}\quad \textrm  {donc} \quad h'(n)<{{a+n-1}\choose {n}}.$$
   Le lemme \ref{lemme technique} appliqué à
$h(n+1)={{a+n-1}\choose {n+1}}+h'(n)$ (avec $i=n+1$ et $j=n$) donne :
 $$h(n+1)^{<n+1>}={{a+n}\choose {n+2}}+h'(n)^{<N>}.$$

Alors puisque $h$ est une fonction de Macaulay on a  :
$$h(n+2)={{a+n}\choose
{n+2}}+h'(n+1)\leq h(n+1)^{<n+1>}={{a+n}\choose {n+2}}+h'(n)^{<N>}$$
d'où le résultat.

\begin{itemize}
 \item
Si $h'(1)<a$, on pose $h'=h_1$ et on a $h=h_0+h_1[-1]$.
\item Si $h'(1)=a$, $
h'^{\sharp}(\infty)=h^{\sharp}(\infty)-h_0^{\sharp}(\infty)
<h^{\sharp}(\infty)$ et on peut appliquer l'hypothèse de récurrence à $h'$ ;
  il
existe des fonctions à support fini  $h_1, \ldots,h_r$ avec $r\geq 2$ et :
 
  \begin{itemize}
\item $h_i(1)=a-1$ pour $0\leq i<r$, $h_r(1)\leq a-1$,
  \item $\sup h_i <s_0(h_{i-1})-1$,
\item $h'=h_1+h_2[-1]+\cdots +h_r[1-r]$.
\end{itemize}
 
  On a alors  $h=h_0+h_1[-1]+\cdots +h_r[-r]$ et
il reste à vérifier que $\sup h_1 <s_0(h_{0})-1$.  Or on a $\sup h_1\leq \sup h'\leq
N-2$ et $s_0(h_{0})=N$ d'où le résultat.
 \end{itemize}

\end{itemize}
 
 \subsubsection*{existence dans le cas général}
  S'il existe $n \in
{\mathbb N }$ avec $ h(n)< {{a+n-2}\choose{n}}$, on définit $N$, $h_0$ et $h'$ comme précédemment et on montre de la même manière que $h_0$ et $h'$ sont des fonctions de Macaulay. De plus, $h'$ étant à support fini, on peut lui appliquer le résultat précédent.
 
 Sinon, l'ensemble :
 $$\bigg\{s \in
{\mathbb N }\mid \forall n, h(n)\geqslant {{a+n-2}\choose{n}}+ {{a+n-3}\choose{n-1}}+\cdots + {{a+n-s-2}\choose{n-s}}\bigg\}$$
 est (un intervalle) non vide et borné supérieurement.  
 
 En effet, sinon il est égal à ${\mathbb N }$, donc $\forall s\in {\mathbb N }$ et $\forall s\in {\mathbb N }$, on a :
 $$h(n)\geqslant {{a+n-2}\choose{n}}+ {{a+n-3}\choose{n-1}}+\cdots + {{a+n-s-2}\choose{n-s}}.$$
 En particulier, pour $s=n$, on obtient :
 $$h(n)\geqslant  {{a+n-2}\choose{n}}+ {{a+n-3}\choose{n-1}}+\cdots + {{a-2}\choose{0}}= {{a+n-1}\choose{n}}$$ donc $h(n)= {{a+n-1}\choose{n}}$ et $h$ est binômiale.
 
 Soit alors $t$ la borne supérieure de cet intervalle. On 
 pose  :
 $$h(n)=\sum_{i=0}^t {{a+n-i-2}\choose{n-i}}+h'(n-t-1)$$
 
 et on montre comme dans la preuve de  \ref {theo principal} que $h'$ est de Macaulay de type $\leqslant a$. 
 
 Si $h'(1)<a$, l'égalité ci-dessus est la décomposition cherchée et vérifie les conditions i), ii) et iii) avec $r=t+1$, $h_i(n)=  {{a+n-i-2}\choose{n-i}}$ pour $0\leqslant i\leqslant t$ et $h_r=h'$.
 
 Si $h'(1)=a$, par définition de $t$, il existe $n \in
{\mathbb N }$ avec $ h'(n)< {{a+n-2}\choose{n}}$ et on peut appliquer 
 \ref {theo principal} à $h'$ ;  il existe $r+1$ fonctions de
Macaulay ($r\geq 1$) $h_0, \ldots,h_r$ avec :
\begin{itemize}
\item[i)] $h_i(1)=a-1$ pour $0\leq i<r$, $h_r(1)\leq a-1$,
\item[ii)] pour $1\leqslant i\leqslant r$, $h_i$ est à support fini et $\sup h_i <s_0(h_{i-1})-1$,
\item[iii)] $h=h_0+h_1[-1]+\cdots +h_r[-r]$.
 \end{itemize}
 
 On a alors $$h(n)=\sum_{i=0}^t {{a+n-i-2}\choose{n-i}}+h_0(n-t-1)+\cdots +h_r(n-t-r)$$
 et on vérifie que les conditions de l'énoncé sont encore vérifiées pour cette décomposition.

 \begin{rem}
\rm

Dans la décomposition obtenue $h=h_0+h_1[-1]+\cdots +h_r[-r]$ on distinguera :
\begin{itemize}
\item si $h(n)$ est d'ordre $a-2$ en $n$ quand $n$ tend vers $+\infty$,
des fonctions de Macaulay  $h_0,\dots h_t$ binômiales de type $a-1$ ;
\item des fonctions de Macaulay $h_{t+1},\dots h_r$ non binômiales, à support fini et de type $a-1$, sauf éventuellement $h_{t+1}$ qui peut être à support infini et $h_r$ qui peut être de type $<a-1$.
\end{itemize}
Si $h$ est à support fini, les fonctions $h_i$ sont toutes à support fini.
 \end{rem}

\end{proof}
 
\section{Sous-schémas ACM de codimension 3}
Nous allons  nous limiter aux sous-schémas ACM de ${\mathbb P}^N$. Cependant  nous pouvons, en appliquant  \ref{proprietes de gamma} aux sous-schémas de codimension 3, énoncer le résultat suivant :
\begin{pro}
\label{proprietes de gamma_X}
Soit $X$ un sous-schéma de pure  codimension 3 de 
${\mathbb P}^N$ et $ \gamma=
\gamma_X$ son caractère de postulation. Il a les propriétés suivantes :
\begin{enumerate}
\item[i)] $\gamma(n) = 0$ pour $n<0$,
\item [ii)]$\gamma(n) = -(n+1)$ pour $0 \leq n < s_0(X) = \inf
\{\, n \in  {\mathbb Z } \mid h^0{\mathcal I}_X(n) \neq 0 \,\}$,
\item [iii)]$\gamma(s_0) \geq -s_0$,
\item[iv)]  $s_1= \inf
\{\, n \geq s_0 \mid \gamma_C(n) >-s_0 \,\}$.
\end{enumerate}
\end{pro}
 
\begin{pro}
Soit $X$ un sous-schéma ACM de 
${\mathbb P}^N$ de codimension 3,
trac\'e sur un sous-schéma  ACM, $Y$, de codimension 2. Alors la fonction  $\gamma_X-\gamma_Y^{\sharp}$ est positive.
\end{pro}

\begin{proof}
On peut supposer, quitte à faire un changement de coordonnées, que le plan défini par $(X_0,\dots ,X_{r-2})$ ne rencontre pas $Y$. Soit $S'= k[X_0,X_1,\dots,X_{r-2}]$. L'anneau $S/{I_X} $ est de profondeur $2$, donc de dimension projective 1, sur $S'$,  et l'idéal $I_{X/Y} $ est un $S'$-module libre gradué. De plus, l'anneau gradué $S/I_Q$ est un $S'$ module libre gradué. On a donc des isomorphismes :
$$S/I_Y\simeq \oplus_{{n \in {\mathbb Z }}}S'(-n)^{L_0(n} \qquad I_{X/Y}\simeq  \oplus_{{n \in {\mathbb Z }}}S'(-n)^{L(n)}$$ 
où $L$ et $L_0$ sont deux fonctions de ${\mathbb N }$ dans  ${\mathbb N}$.

On en déduit :
  \begin{equation*}
  \begin{split}
  h^0{\mathcal O}_{X}(n)&=\sum_{m \in {\mathbb Z }}{{n-m+N-2}\choose {N-2}}(L_0(m)-L(m))\\
 \gamma_X(n)&=-\partial^{N-1}h^0{\mathcal O}_{X}(n)=-\sum_{m \in {\mathbb Z }}{{n-m-1}\choose -1}(L_0(m)-L(m))\\
  \gamma_X&=L-L_0\\
  \end{split}
    \end{equation*}
   
et :
 \begin{equation*}
 \begin{split}
h^0{\mathcal O}_{Y}(n)&=\sum_{m \in {\mathbb Z }}{{n-m+N-2}\choose {N-2}}L_0(m)\\
 \gamma_Y^{\sharp}(n)&=-\partial^{N-1}h^0{\mathcal O}_{Y}(n)=-\sum_{m \in {\mathbb Z }}{{n-m-1}\choose -1}L_0(m)=-L_0(n)\\
 \end{split}
   \end{equation*}
   
   donc $\gamma_X-\gamma_Y^{\sharp}=L$ est positive.

\end{proof}
\begin{rem}
\rm Ce résultat est  déjà vrai si $X$ est de codimension 2 et $Y$ de codimension 1.  Dans ce cas si $d$ est le degré de $Y$, la fonction $\gamma_Y^{\sharp}$ vaut $-1$ sur l'intervalle $[0,d-1]$ et 0 ailleurs. La fonction $\gamma_X$ vaut $-1$ sur l'intervalle $[0,s_0(X)-1]$ et est positive pour $n\geq s_0(X)\geq d$, d'où le résultat.
\end{rem}

% courbes integres

\begin{coro}
\label{ACM integres}
Soit $X$ un sous-schéma ACM intègre de 
${\mathbb P}^N$ de codimension 3, $\gamma_X$ son caractère, $s_0=s_0(X)$ et $s_1= \inf
\{\, n \geq s_0 \mid \gamma_X(n) >-s_0 \,\}$. Alors on a, pour $n\geq s_1$, $\gamma_X(n)\geq \inf \{0, n-s_0-s_1+1\}$.

\end{coro}

\begin{proof}
Soit $s_1= \inf
\{\, n \geq s_0 \mid \gamma_X(n) >-s_0 \,\}$.

Soit $Y$ l'intersection complète  de deux hypersurfaces de degrés $s_0$ et $s_1$.
On calcule le caractère $\gamma_Y$ grâce à la résolution : 
$$0\to S[-s_0-s_1] \to S[-s_0]\oplus S[-s_1] \to I_Y \to 0$$ et on obtient :
\begin{equation*}
\gamma_Y(n) =
\begin{cases}& -1  \textrm { pour } 0\leq n<s_0\\
& 1  \textrm { pour } s_1\leq n<s_0+s_1\\
&  0 \textrm { sinon }
\end{cases}
\end{equation*}

\begin{equation*}
\gamma_Y^{\sharp}(n) =
\begin{cases}& -(n-1)  \textrm { pour } 0\leq n<s_0\\
& -s_0  \textrm { pour } s_0\leq n<s_1\\
& n-s_0-s_1+1  \textrm { pour } s_1\leq n<s_0+s_1-1\\
& 0 \textrm { sinon }
\end{cases}
\end{equation*}

Donc $\gamma_X-\gamma_Y^{\sharp}$ est toujours nulle pour $n<s_1$. 

Si $X$ est intègre, il existe  une hypersurface intègre $W$ de degré $s_0$ contenant $X$, et une hypersurface $W'$ contenant $X$ de degré $s_1$ ne contenant pas $W$. 
Alors $X$ est contenue dans l'intersection complète $Y$ de $W$ et $W'$ et $\gamma_X-\gamma_Y^{\sharp}$ est positive.

On en déduit pour $n\geq s_1$, $\gamma_X(n)\geqslant-\gamma_Y^{\sharp}(n)=\inf \{0, n-s_0-s_1+1\}$.

\end{proof}
% caractères de postulation des courbes ACM de ${\mathbb P}^4$.
Nous pouvons, grâce à \ref{theo principal}, décrire tous les caractères de postulation des sous-schémas ACM  de 
${\mathbb P}^N$ de codimension 3.

\begin{pro}
\label{courbesACM}
Soit $X$ un sous-schéma ACM de 
${\mathbb P}^N$ de codimension 3 et $\gamma_X $ son caractère de postulation.
Il existe un entier $r\geqslant 0$, des caractères  positifs $\gamma_0$, $\gamma_1$, \ldots,
$\gamma_r$ (cf. \ref{definition caractere positif}) vérifiant :
 \begin{itemize}
\item $r=s_0(X)-1$,
 \item
$\sup \gamma_i <s_0(\gamma_{i-1})$,
 \item $\gamma_X=\gamma_0+\gamma_1[-1]+\cdots +\gamma_r[-r]$. 
 \end{itemize}
 Inversement, pour tout entier $r\leqslant0$, pour tous
caractères positifs $\gamma_0$, $\gamma_1$, \ldots,
$\gamma_r$ vérifiant 
$\sup \gamma_i <s_0(\gamma_{i-1})$, il existe un sous-schéma ACM  de 
${\mathbb P}^N$ de codimension 3, $X$,
tel que $\gamma_X=\gamma_0+\gamma_1[-1]+\cdots +\gamma_r[-r]$.
 \end{pro}
\begin{proof}
Si $X$ est contenu dans un hyperplan, son caractère est celui d'un sous-schéma ACM  de 
${\mathbb P}^{N-1}$ de codimension 2. Il est positif. On pose $\gamma_X=\gamma_0$.

  Si $X$ n'est pas dégénéré, son h-vecteur $h_X$ est de type 3. On lui applique le théorème \ref{theo principal}  : il existe 
 $r+1$ fonctions de
Macaulay à support fini ($r\geq 1$) $h_0, \ldots,h_r$ avec :
\begin{itemize}
\item $h_i(1)=2$ pour $0\leq i<r$, $h_r(1)\leq 2$,
\item $\sup h_i <s
  _0(h_{i-1})-1$,
\item $h_X=h_0+h_1[-1]+\cdots +h_r[-r]$.
\end{itemize}
D'après \ref{caractere positif-fonction de Macaulay} pour tout $i$, $\gamma_i=-\partial h_i$ est un caractère
positif,  $\sup \gamma_i=\sup h_i+1$ et $s_0(\gamma_i)=s_0(h_i)$ pour $i<r$.

La propriété inverse est la conséquence de \ref{Migliore-Nagel}.
 
 \end{proof}
 
 \begin{pro}
Avec les notations de \ref{courbesACM}, on a les propriétés suivantes :
\begin{itemize}
\item pour $s_0(X)\leq n< s_0(\gamma_{r-1})+s_0(X)-2$,  on a  $\gamma_X(n)\geqslant-s_0(X)$,
\item pour $s_0(\gamma_{i})+i\leq n <s_0(\gamma_{i-1})+i-1$, on a 
$\gamma_X(n)\geqslant-i$,
\item pour $s_0(\gamma_{0})\leq n$, on a $\gamma_X(n)\geqslant0$.

\end{itemize}
\end{pro}

\begin{proof}
Si $X$ est dégénéré, on a $r=0$, $\gamma_X=\gamma_0$ est donc un caractère positif.

Supposons $X$ non dégénéré et gardons les notations de la démonstration de \ref{courbesACM}. Posons $s_0=s_0(X)(=r+1)$. On a vu dans la preuve de  \ref{proprietes de h} que, pour
 $s_0\leq n< s_0(h_{r-1})+r-1$, on a :
$$h_X(n)={{n+1}\choose{1}}+\cdots
+{{n-r+2}\choose{1}}+h_r(n-r)$$
On en déduit qu'on a : \begin{equation*}
\begin{split}
\gamma_X(s_0)=\gamma_X(r+1)&={{r+2}\choose{2}}-{{r+2}\choose{1}}-\cdots
-{{3}\choose{1}}-h_r(1)\\
&={{r+2}\choose{2}}-{{r+3}\choose{2}}+3-h_r(1)\\
&=1-r-h_r(1)=\gamma_r(1)-r.\\
\end{split}
\end{equation*}

De même, pour  $s_0<n< s_0(h_{r-1})+r-1$ on a :
$$\gamma_X(n)=-{{n}\choose{0}}-\cdots
-{{n-r+1}\choose{0}}+\gamma_r(n-r)=\gamma_r(n-r)-r$$
et cette formule est donc valable pour $n=s_0$.

De même, pour $s_0(\gamma_{i})+i\leq n <s_0(\gamma_{i-1})+i-1$, on a :
$$h_X(n)={{n+1}\choose{1}}+\cdots
+{{n+2-i}\choose{1}}+h_i(n-i)$$
donc pour $s_0(\gamma_{i})+i<n <s_0(\gamma_{i-1})+i-1$, on a :
   \begin{equation*}
\begin{split}
\gamma_X(n) &=-{{n}\choose{0}}-\cdots
-{{n-i+2}\choose{0}}+\gamma_i(n-i)\\
&=\gamma_i(n-i)-1\geqslant -i.\\
\end{split}
   \end{equation*}
   
 On laisse le soin au lecteur de vérifier que c'est encore vrai pour $n=s_0(\gamma_{i})+i$.
 
 Pour $n>s_0(\gamma_{0})$, on a $\gamma_X(n)=\gamma_0(n)\geqslant 0$ et on vérifie que c'est encore vrai pour $n=s_0(\gamma_{0})$.
 
 \end{proof}

\begin{coro}
Avec les notations de \ref{courbesACM}, on a $$s_1(X)= \inf
\{\, n \geq s_0 \mid \gamma_X(n) >-s_0 \,\}=s_0(\gamma_{r})+s_0(X)-1.$$
\end{coro}

\begin{proof}
Si l'intervalle $ [r+1,s_0(h_{r-1})+r-1[ $ n'est pas vide,
on a en fait $s_1(X)= \inf
\{\,  s_0(X)\leqslant n\leqslant s_0(h_{r-1})+r-1\mid \gamma_X(n) >-s_0 \,\}$. Sur cet intervalle on a 
$ \gamma_X(n)= \gamma_r(n-r)-s+1$, donc $s_1(X)= s_0(\gamma_{r})+s_0(X)-1$.

Si  l'intervalle $ [r+1,s_0(h_{r-1})+r-1[ $ est vide, pour tout $n>s_0=r+1$ on a $\gamma_X(n) >-s_0$. De plus dans ce cas $s_0(h_{r-1})=2$, $\sup (h_r)=0$, et $s_0(\gamma_{r})=1$. La formule est donc encore valable.
\end{proof}

\begin{rem}
\label{astuce}
\rm
Dans la pratique,  on  construit les $\gamma_i$ de proche en proche, comme on l'a fait dans la démonstration de \ref{theo principal}. Soit 
$N=\inf \{\, n \in
{\mathbb N }
\mid \sum_{m>n}\gamma (n)<n \,\}$.

On définit  $\gamma_0$ par :
 
 \begin{equation*}
\gamma_0(n) =
\begin{cases} -1 & \textrm{pour } n<N\\
\gamma(n) & \textrm{pour }n> N\\
 \end{cases}\qquad
  \gamma_0(N)=\gamma(N)- \sum_{m>N}\gamma (n) 
\end{equation*}

et on recommence.
\end{rem}

  \begin{rem}
  \label{existence}
\rm Soient  $\gamma_0$, $\gamma_1$, \ldots,
$\gamma_r$ des caractères  positifs vérifiant les conditions de l'énoncé de   \ref{courbesACM}.
  Supposons qu'on puisse construire une suite de sous-schémas ACM de codimension 3, $X_i$, ($0\leqslant i\leqslant r$) de la manière suivante :
  \begin{itemize}
  
 \item  $X_r$ est un sous-schéma hyperplan de codimension 3 de caractère $\gamma_r$,
 \item pour $0\leqslant i <r$, $X_i$ est obtenu à partir de $X_{i+1}$ par une
biliaison \'el\'ementaire Gorenstein de hauteur $1$ sur un sous-schéma  ACM de codimension 2 $Y_i$ de caractère  $\gamma_i$ (cf.  \ref{bil-elem}).
\end{itemize}

Alors d'après \ref{variation Gorenstein}, le caractère de $X_0$ est  $\gamma_{X_0}=\gamma_0+\gamma_1[-1]+\cdots +\gamma_r[-r]$.\\
Si on sait construire une telle suite, on redémontre dans ce cas le résultat d'existence  \ref{Migliore-Nagel}.
Le problème est à chaque pas de trouver un sous-schéma  ACM de codimension 2  de caractère  $\gamma_i$ contenant $X_{i+1}$.

\end{rem}

%Courbes ACM tracées sur une quadrique

  \subsection*{Sous-schémas ACM de codimension 3 tracés sur une hypersurface quadrique}

Nous allons décrire tous les caractères de  postulation des sous-schémas ACM de codimension 3, $X$, vérifiant $s_0(X)=2$.

\begin{pro}
Soit $X$ un sous-schéma ACM de codimension 3 non dégénéré de ${\mathbb P}^N$
tracé sur une quadrique et  $\gamma_X=\gamma$ son caractère. Les propriétés équivalentes suivantes  sont réalisées :

\begin{itemize}
\item {i)} il existe deux caractères  positifs (cf. \ref{definition caractere positif}) $\gamma_0$ et $\gamma_1$ vérifiant 
$\sup \gamma_1<s_0(\gamma_{0})$, tels qu'on ait 
$\gamma_X=\gamma_0+\gamma_1[-1]$,

\item {ii)} il existe deux entiers $1\leqslant t<s$ tels qu'on ait :

\begin{equation*}
\gamma(n)
\begin{cases}
      &=-2 \text{ pour } 1\leqslant n\leqslant t, \\
      &\geqslant -1\textrm { pour }  t<n< s, \\
      & \geqslant 0\text{ pour }n\geqslant s
\end{cases}
\quad \rm{et}\quad \sum_{n>s}\gamma (n)\leqslant s \leqslant  \sum_{n\geqslant s}\gamma (n)
\end{equation*}.
\end{itemize}

Inversement tout caractère $\gamma$ vérifiant i) ou ii)  est le caractère d'un sous-schéma ACM de codimension 3  non dégénéré $X$ de ${\mathbb P}^N$
tracé sur une quadrique. De plus on a $s_1(X)=s_0(\gamma_1)+1=t+1$.
\end{pro}

\begin{proof}

i) est une conséquence de \ref{courbesACM}. 

Posons $s=s_0(\gamma_0)$ et $t=s_0(\gamma_1)$. On a $t\leqslant \sup \gamma_1<s$.

Pour $1\leqslant n\leqslant t$, $\gamma(n)=-2$.

Pour $t+1\leqslant n<s$, $\gamma(n)=-1+\gamma_1(n-1)\geqslant -1$.

Pour $n=s$, $\gamma(s)=\gamma_0(s)+\gamma_1(s-1)\geqslant 0$.

Pour $s< n$, $\gamma(n)=\gamma_0(n)\geqslant 0$.

De plus, on a $$ \sum_{n>s}\gamma (n)= \sum_{n>s}\gamma_0 (n)\leqslant \sum_{n\geqslant s}\gamma_0 (n)=s$$ et
$$ \sum_{n\geqslant s}\gamma (n)= \sum_{n\geqslant s}\gamma_0 (n)+\gamma_1(s-1)\geqslant \sum_{n\geqslant s}\gamma_0 (n)=s.$$

Inversement, soit $\gamma$ un caractère vérifiant ces propriétés. En s'inspirant  de \ref{astuce}, on définit  $\gamma_0$ et $\gamma_1$ de la manière suivante :

\begin{equation*}
\gamma_0(n)=\begin{cases}
      & -1 \textrm{ pour } n<s, \\
      & \,\gamma(n) \textrm{ pour }n>s\\
      & s-\sum_{n>s}\gamma (n)\textrm{ pour } n=s
\end{cases}
\quad \rm{et}\quad \gamma_1=\gamma[1]+\gamma_0[1].
\end{equation*}

On vérifie que $\gamma_0$ et $\gamma_1$ sont des  caractères positifs, et qu'on a $\sup \gamma_1<s_0(\gamma_{0})$.

\end{proof}

En utilisant \ref{ACM integres}, on peut décrire les caractères des sous-schémas intègres ACM de codimension 3  de ${\mathbb P}^N$
tracés sur une quadrique :
\begin{coro}
\label{coro integre}
Soit $X$ un sous-schéma intègre ACM de codimension 3  de ${\mathbb P}^N$
tracé sur une quadrique et  $\gamma_X=\gamma$ son caractère. Il existe un entier $1\leqslant t$ tel qu'on ait :

\begin{equation*}
\gamma(n)
\begin{cases}
      &=-2 \text{ pour } 1\leqslant n\leqslant t, \\
      &\geqslant -1\textrm { pour } n=t+1, \\
      & \geqslant 0\text{ pour }n\geqslant t+2.
\end{cases}
\end{equation*}

\end{coro}
\begin{proof}
 On sait que pour $n\geq s_1(X)$, $\gamma_X(n)\geq \inf \{0, n-s_0(X)-s_1(X)+n+1\}$ d'où le résultat puisque $s_0(X)=2$ et $s_1(X)=t+1$.
 \end{proof}

 En fait, dans le cas des sous-schémas ACM de codimension 3  de ${\mathbb P}^N$ tracés sur une quadrique, on sait redémontrer directement le résultat de \ref{Migliore-Nagel}. Plus précisément : 
 
 \begin{pro}
 \label{surfaces emboitees}
 Soient $\gamma_0$ et $\gamma_1$ deux caractères  positifs vérifiant 
$$\sup \gamma_1\leqslant  \inf
\{\, n\in \mathbb{N} \mid \gamma_0(n) >0 \,\}.$$   Il existe un sous-schéma  ACM de codimension 2,  $Y_1$, de caractère $\gamma_1$ et un sous-schéma  ACM de codimension 2, $Y_0$, de caractère $\gamma_0$ qui le contient.
 
\end{pro}

\begin{proof}
La démonstration se fait par récurrence sur l'entier $\rm {deg}\,\gamma_0+\rm {deg}\,\gamma_1$.
\begin{itemize}
\item Si $\rm {deg}\,\gamma_0+\rm {deg}\,\gamma_1=2$, on choisit pour $Y_0$ et $y_1$ le même $(N-2)$-plan $P$.
\item Si $\rm {deg}\,\gamma_0+\rm {deg}\,\gamma_1>2$, on distingue deux cas :
 \begin{itemize}
\item si $\sup \gamma_1< \inf
\{\, n\in \mathbb{N} \mid \gamma_0(n) >0 \,\}=s$, on peut faire à partir de tout sous-schéma  ACM de codimension 2  de caractère $\gamma_0$ une biliaison élémentaire (intersection complète) descendante de hauteur $-1$ sur une hypersurface de degré $s$ \cite{MD}. Soit $\gamma'_0$ le caractère obtenu. Alors on a $\rm{deg}\,\gamma'_0 <\rm {deg}\,\gamma_0$ et 
:
$$ \sup \gamma'_1\leqslant s-1\leqslant  \inf
\{\, n\in \mathbb{N} \mid \gamma'_0(n) >0 \,\}.$$
D'après l'hypothèse de récurrence, il existe un sous-schéma  ACM de codimension 2, $Y_1$, de caractère $\gamma_1$ et un sous-schéma  ACM de codimension 2, $Y'_0$, de caractère $\gamma'_0$ qui le contient. Si $W$ est une  hypersurface de degré $s$ contenant $Y'_0$, on choisit pour $Y_0$ la réunion de $Y'_0$ et d'une section hyperplane de $W$. C'est un sous-schéma   obtenu par une biliaison élémentaire (intersection complète) triviale de hauteur $1$ à partir de $Y'_0$ et il a pour caractère $\gamma_0$.

 \item si $ \sup \gamma_1= \inf
\{\, n\in \mathbb{N} \mid \gamma_0(n) >0 \,\}=s$, soient comme précédemment $\gamma'_0$ et $\gamma'_1$ les caractères obtenus par une biliaison élémentaire (intersection complète) descendante de hauteur $-1$ sur une hypersurface de degré $s$ à partir de $\gamma_0$ et $\gamma_1$. On a $\rm{deg}\,\gamma'_0 <\rm {deg}\,\gamma_0$, $\rm{deg}\,\gamma'_1 <\rm {deg}\,\gamma_1$ et : 
 $$  \sup
\{\, n\in \mathbb{N} \mid \gamma'_1(n) >0 \,\}\leqslant s-1\leqslant  \inf
\{\, n\in \mathbb{N} \mid \gamma'_0(n) >0 \,\}.$$
 D'après l'hypothèse de récurrence, il existe un sous-schéma  ACM de codimension 2, $Y'_1$, de caractère $\gamma'_1$ et un sous-schéma  ACM de codimension 2, $Y'_0$, de caractère $\gamma'_0$ qui le contient. Soit $W$ une  hypersurface de degré $s$ contenant $Y'_0$ (et $Y'_1$). On choisit pour $Y_0$ (resp. $Y_1$) la réunion de $Y'_0$ (resp. $Y'_1$) et d'une section hyperplane de $W$. 
 \end{itemize}

\end{itemize}
\end{proof}

\begin{coro}
 Soient $\gamma_0$ et $\gamma_1$ deux caractères  positifs vérifiant 
$\sup \gamma_1<s_0(\gamma_{0})$.   Il existe un sous-schéma  ACM de codimension 3, $X$, de caractère $\gamma_X=\gamma_0+\gamma_1[-1]$ obtenu par biliaison \'el\'ementaire Gorenstein de hauteur $1$ sur un sous-schéma  ACM de codimension 2, $Y_0$, de caractère  $\gamma_0 $  à partir d'un sous-schéma  ACM de codimension 3 hyperplan de caractère $\gamma_1$.
\end{coro}

\begin{proof}
Les conditions de  \ref{surfaces emboitees} étant satisfaites, il existe un sous-schéma  ACM de codimension 2, $Y_1$, de caractère $\gamma_1$ et un sous-schéma  ACM de codimension 2, $Y_0$, de caractère $\gamma_0$ qui le contient. Soit $X_1$ une section hyperplane de $Y_1$. On peut faire une biliaison  \'el\'ementaire Gorenstein de hauteur $1$ sur $Y_0$ à partir de $X_1$. D'après  \ref{existence} le sous-schéma obtenu a pour caractère $\gamma_0+\gamma_1[-1]$.
\end{proof}

\begin{exa}
\rm Dans le cas où $\gamma_1$ est le caractère d'une droite ($\gamma_1=(-1,1)$), on peut donner une autre description d'une courbe  de  ${\mathbb P}^4$ de caractère $\gamma_0+\gamma_1[-1]$ de la manière suivante : soient $C'$ une courbe ACM contenue dans un hyperplan $H$ de  ${\mathbb P}^4$ de caractère $\gamma_0$ et $D$ une droite non contenue dans $H$ et coupant $C'$ en un point. La réunion $C$ de $C'$ et $D$ est une courbe ACM tracée sur une quadrique. De la  suite exacte :
$$0\to \mathcal{I}_{C}\to \mathcal{I}_{C'}\to \mathcal{O}_{D}(-1) \to 0$$
on déduit :
$$h^0\mathcal{I}_{C}(n)=h^0\mathcal{I}_{C'}(n)-h^0\mathcal{O}_{D}(n-1)$$
$$\gamma_{C}(n)=\gamma_{C'}(n)-\partial^3h^0\mathcal{O}_{D}(n-1)=\gamma_{C'}(n)+\gamma_{D}(n-1).$$

D'après \ref{coro integre}, si $s_0(C')>2$, aucune courbe intègre de ${\mathbb P}^4$ ne peut avoir ce caractère. 
\end{exa}

\begin{exa}
\rm Soient $P_1$ et $P_2$ deux plans de ${\mathbb P}^4$ se coupant en un point $Z$, $C_1$ et $C_2$ deux courbes planes contenues respectivement dans $P_1$ et $P_2$ et se coupant en un point. La réunion $C$ de $C_1$ et $C_2$ est une courbe ACM tracée sur une quadrique.

Des  suites exactes :
$$0\to \mathcal{I}_{C}\to \mathcal{I}_{C_1}\to \mathcal{I}_{Z/{C_2}} \to 0$$
$$ 0\to \mathcal{I}_{Z/{C_2}} \to \mathcal{O}_{C_2} \to \mathcal{O}_{Z}\to 0$$
on déduit :
$$h^0\mathcal{I}_{C}(n)=h^0\mathcal{I}_{C_1}(n)+h^0\mathcal{I}_{C_2}(n)-h^0\mathcal{O}_{{\mathbb P}}(n)+1$$
$$\gamma_{C}=\gamma_{C_1}+\gamma_{C_2}+(1,-2,1).$$

De plus on rappelle (cf. \ref{ex.hypersurface}) que si $d_1$ est le degré de $C_1$, $\gamma_{C_1}$ ne prend que deux valeurs non nulles, $\gamma_{C_1}(0)=-1$, $\gamma_{C_1}(d_1)=1$, et il en est de même de  $\gamma_{C_2}$.

On remarque que si $D$ est une droite, $\gamma_0=\gamma_{C}-\gamma_{D}[-1]$ est un caractère positif, donc que $C$ a le même caractère que la réunion d'une courbe hyperplane de caractère $\gamma_0$ et d'une droite la coupant en un point et ``sortant'' de l'hyperplan de la courbe.

\end{exa}

%petits (d,g)
\section {Courbes ACM de ${\mathbb P}^4$}

\begin{pro}
Avec les notations de \ref{courbesACM}, désignons par $Y_i$ une surface  ACM de 
caractère $\gamma_i$. On a :
$$d=\sum_{i=0}^r d_i \quad \quad 2g-2=\sum_{i=0}^r(\delta_i+(2i+1)d_i).$$
\end{pro}
 \begin{proof}
  On rappelle qu'on a 
:
   
$$d=\sum _{k \in {\mathbb Z }}k\gamma_C(k) \quad \quad g-1=\sum _{k \in {\mathbb Z }}
{{(k-1)(k-2)\over 2}}\gamma_C(k).$$
  
  Alors 
   \begin{equation*}
 \begin{split}
d&=\sum _{k ,i}k\gamma_i(k-i)=\sum _{k ,i}(k+i)\gamma_i(k)\\
&=\sum_{k,i} k\gamma_i(k)+\sum_{k,i} i\gamma_i(k)=\sum_i d_i\\
  \end{split}
  \end{equation*}
 
puisque $\gamma_i$ est un caractère ($\sum_{k} k\gamma_i(k)=0$).
  
   \begin{equation*}
 \begin{split}
 2g-2&=\sum _{k,i} {(k-1)(k-2)}\gamma_i(k-i)=\sum _{k,i} {(k+i-1)(k+-2)}\gamma_i(k)\\
 &=\sum _{k,i} k^2\gamma_i(k)+\sum _{k,i} (2i-3)k\gamma_i(k)+\sum _{k,i} (i-1)(i-2)\gamma_i(k)\\
 &=\sum_i (\delta_i+4d_i)+\sum_i (2i-3)d_i=\sum_{i}(\delta_i+(2i+1)d_i).\\
  \end{split}
   \end{equation*}

  \end{proof}

  \subsection*{Courbes ACM de degré $\leqslant 10$}

Pour terminer, nous allons utiliser \ref{courbesACM} pour donner la liste des degrés et genres des courbes ACM  de ${\mathbb P}^4$ de degré  $\leqslant 10$, ce qui complète la liste donnée par \cite {H1}. Nous préciserons à chaque fois dans la preuve les caractères $\gamma_i$ et $\gamma_C$.

\begin{pro}
Soit $C$ une courbe ACM  non dégénérée de ${\mathbb P}^4$ de degré  $\leqslant 10$, $d$ son degré et $g$ son genre. Alors le couple $(d,g)$ prend ses valeurs dans l'ensemble suivant : \\
$(4,0)$, $(5,1)$, $(6,2)$, $(6,3)$,  $(7,3)$, $(7,4)$, $(7,6)$, $(8,4)$, $(8,5)$, $(8,6)$, $(8,7)$, $(8,10)$, $(9,5)$, $(9,6)$, $(9,7)$, $(9,8)$, $(9,9)$, $(9,11)$, $(9,15)$, $(10,6)$, $(10,7)$,$(10,8)$, $(10,9)$, $(10,10)$, $(10,12)$,  $(10,13)$ et  $(10,16)$.  
\end{pro}

\begin{proof}

Grâce à la borne $d\geqslant {{s+2}\choose {3}}$ (cf. \ref{borne inf}), on a  $2\leqslant s_0(C)\leqslant 3$.

\begin{itemize}
\item  $s_0(C)=3$. Il y a un seul caractère de degré  $\leqslant 10$, qui réalise la borne inférieure de \ref{borne inf}. On a :
$\gamma_0=(-1,-1,-1,3)$,
$\gamma_1=(-1,-1,2)$
 et $\gamma_C=(-1,-2,-3,6)$, $(d,g)=(10,6)$.

\item $s_0(C)=2$.
\begin{enumerate}
\item $\gamma_1=(-1,0,1)$. On a  les possibilités suivantes :
\begin{enumerate}
\item  $\gamma_0=(-1,-1,-1,3)$, $\gamma_C=(-1,-2,-1,4)$ et $(d,g)=(8,4)$ ;
\item $\gamma_0=(-1,-1,-1,2,1)$, $\gamma_C=(-1,-2,-1,3,1)$ et $(d,g)=(9,6)$ ;
\item $\gamma_0=(-1,-1,-1,2,0,1)$, $\gamma_C=(-1,-2,-1,3,0,1)$ et $(d,g)=(10,9)$ ;
\item $\gamma_0=(-1,-1,-1,1,1,1)$, $\gamma_C=(-1,-2,-1,2,2)$ et $(d,g)=(10,8)$.
\end{enumerate}

\item $\gamma_1=(-1,-1,2)$. On a  les possibilités suivantes :
\begin{enumerate}
\item $\gamma_0=(-1,-1,-1,3)$, $\gamma_C=(-1,-2,-2,5)$ et $(d,g)=(9,5)$ ;
\item $\gamma_0=(-1,-1,-1,2,1)$, $\gamma_C=(-1,-2,-2,4,1)$ et $(d,g)=(10,7)$.

\end{enumerate}
\item $\gamma_1=(-1,1)$.
On a alors nécessairement 
 $\gamma_0=(-1,-1)+1_{[a]}+1_{[b]}+1_{[c]}$ avec $2\leqslant a\leqslant b\leqslant c$ et $a+b+c\leqslant 12$, donc  $\gamma_C=(-1,-2)+1_{[a]}+1_{[b]}+1_{[c]}$, $d=a+b$, $2g=(a-1)(a-2)+(b-1)(b-2)+(c-1)(c-2)$. Lorsque $a,$ $b$, $c$ varient on obtient toutes les autres valeurs annoncées.

\end{enumerate}
\end{itemize}

\end{proof}

\begin{rem}
\rm Il existe une borne supérieure de Castelnuovo pour le genre d'une courbe intègre de ${\mathbb P}^4$, qui montre que certains des caractères précédents ne peuvent pas être réalisés par une courbe intègre, bien que la conditon de  \ref{ACM integres} soit toujours vérifiée, ce qui prouve que ce n'est pas une condition suffisante.
\end{rem}

 {}

\end{document}